\documentclass[12pt,a4paper]{article}
\usepackage{}
\usepackage{indentfirst}
\usepackage{amsfonts}
\usepackage{amssymb}
\usepackage{mathrsfs}
\usepackage{amsmath}
\usepackage{amsthm}
\usepackage{enumerate}
\usepackage{cite}
\usepackage[all]{xy}
\allowdisplaybreaks
\newtheorem{defn}{Definition}[section]
\newtheorem{thm}[defn]{Theorem}
\newtheorem{lem}[defn]{Lemma}
\newtheorem{prop}[defn]{Proposition}
\newtheorem{cor}[defn]{Corollary}
\newtheorem{eg}[defn]{Example}
\newtheorem{re}[defn]{Remark}
\newcommand{\bdefn}{\begin{defn}}
\newcommand{\edefn}{\end{defn}}
\newcommand{\bthm}{\begin{thm}}
\newcommand{\ethm}{\end{thm}}
\newcommand{\blem}{\begin{lem}}
\newcommand{\elem}{\end{lem}}
\newcommand{\bprop}{\begin{prop}}
\newcommand{\eprop}{\end{prop}}
\newcommand{\bcor}{\begin{cor}}
\newcommand{\ecor}{\end{cor}}
\newcommand{\beg}{\begin{eg}}
\newcommand{\eeg}{\end{eg}}
\newcommand{\bre}{\begin{re}}
\newcommand{\ere}{\end{re}}
\newcommand{\bpf}{\begin{proof}}
\newcommand{\epf}{\end{proof}}
\newcommand{\F}{{\rm{\bf F}}}
\newcommand{\id}{{\rm id}}
\newcommand{\Ker}{\rm Ker}

\newcommand{\supercite}[1]{\textsuperscript{\cite{#1}}}
\newcommand{\benu}{\begin{enumerate}}
\newcommand{\eenu}{\end{enumerate}}
\newcommand{\bc}{\begin{center}}
\newcommand{\ec}{\end{center}}
\newcommand{\bea}{\begin{eqnarray}}
\newcommand{\eea}{\end{eqnarray}}
\newcommand{\Bea}{\begin{eqnarray*}}
\newcommand{\Eea}{\end{eqnarray*}}
\newcommand{\beq}{\begin{equation}}
\newcommand{\eeq}{\end{equation}}
\newcommand{\Beq}{\begin{equation*}}
\newcommand{\Eeq}{\end{equation*}}
\newcommand{\bspl}{\begin{split}}
\newcommand{\espl}{\end{split}}
\newcommand\relphantom[1]{\mathrel{\phantom{#1}}}

\begin{document}
\title{\bf  Cohomology and 1-parameter formal deformations of Hom-Lie triple systems}
\author{\normalsize \bf Yao Ma$^1$, Liangyun Chen$^1$,  Jie Lin$^2$}
\date{\small{ $^1$ School of Mathematics and Statistics, Northeast Normal University, Changchun,  130024,  CHINA
  \\$^2$ Sino-European Institute of Aviation Engineering,Civil Aviation University of China,  Tianjin, 300300,  CHINA}}
\maketitle
\begin{abstract}
In this paper, we study the cohomology theory of Hom-Lie triple systems generalizing the Yamaguti cohomology theory of Lie triple systems. We introduce the central extension theory for Hom-Lie triple systems and show that there is a one-to-one correspondence between equivalent classes of central extensions of Hom-Lie triple systems and the third cohomology group. We develop the 1-parameter formal deformation theory of Hom-Lie triple systems, and prove that it is governed by the cohomology group.
\bigskip

\noindent{Key words:}  Hom-Lie triple system, cohomology, central extension, deformation\\
\noindent{Mathematics Subject Classification(2010):}  17A40, 17B10, 17B56, 55U15
\end{abstract}

\footnote[0]{Corresponding author(L. Chen): chenly640@nenu.edu.cn.}
\footnote[0]{Supported by  NNSF of China (No.11171055),  NSF of  Jilin province (No.201115006), Scientific Research Foundation for Returned Scholars Ministry of Education of China  and the Fundamental Research Funds for the Central Universities (No.12SSXT139). }

\section{Introduction}
A Hom-Lie algebra is a vector space endowed with a skew symmetric bracket satisfying a Jacobi identity twisted by a map. Before Hom-Lie algebras appeared, Hu studied $q$-Lie algebras, which are special Hom-Lie algebras\supercite{Hu}. Lie algebras are special cases of Hom-Lie algebras when the twisted map is the identity map. The notion of Hom-Lie algebras was introduced by Hartwig, Larsson and Silvestrov to describe the $q$-deformation of the Witt and the Virasoro algebras\supercite{Hartwig&Larsson&Silvestrov}. Since then, Hom-type algebras have been investigated by many authors\supercite{Ammar&Ejbehi&Makhlouf, Ammar&Mabrouk&Makhlouf, Benayadi&Makhlouf, Elhamdadi&Makhlouf, Larsson&Silvestrov, Liu&Chen&Ma, Makhlouf&Silvestrov, Sheng, Sheng&Chen, Yau1}. In particular, the notion of Hom-Lie triple systems was introduced by Yau\supercite{Yau2}.

A deformation is a tool to study a mathematical object by deforming it into a family of the same kind of objects depending on a certain parameter. The deformation theory was introduced by Gerstenhaber for rings and algebras\supercite{Gerstenhaber1, Gerstenhaber2, Gerstenhaber3, Gerstenhaber4}, and by Kubo and Taniguchi for Lie triple systems\supercite{Kubo&Taniguchi}. They studied 1-parameter formal deformations and established the connection between the cohomology groups and infinitesimal deformations: the suitable cohomology groups for the deformation theory of associative algebras and Lie triple systems are the Hochschild cohomology\supercite{Hochschild} and the Yamaguti cohomology\supercite{Yamaguti}, respectively.

The purpose of this paper is to consider the cohomology theory and the 1-parameter formal deformation theory of Hom-Lie triple systems based on some work in \cite{Kubo&Taniguchi, Sheng, Yamaguti}. The paper is organized as follows. Section 2 is devoted to the cohomology theory of multiplicative Hom-Lie triple systems. Section 3 concerns central extensions. It is proved that there is a one-to-one correspondence between equivalent classes of central extensions of a Hom-Lie triple system and the third cohomology group with coefficients in the trivial representation. Section 4 is dedicated to the 1-parameter formal deformation theory of multiplicative Hom-Lie triple systems. We show that the cohomology group defined in Section 2 is suitable for this 1-parameter formal deformation theory.

Throughout this paper $\F$ denotes an arbitrary field.

\section{The cohomology of Hom-Lie triple systems}
We start by recalling the definitions of Lie triple systems and Hom-Lie triple systems.
\bdefn{\rm\supercite{Lister}}
A vector space $T$ together with a trilinear map $(x, y, z)\mapsto[xyz]$ is called a Lie triple system (LTS) if
\begin{enumerate}[(1)]
\item $[xxz]=0$,
\item $[xyz]+[yzx]+[zxy]=0$,
\item $[uv[xyz]]=[[uvx]yz]+[x[uvy]z]+[xy[uvz]]$,
\end{enumerate}
for all $x,y,z,u,v\in T$.
\edefn

\bdefn{\rm\supercite{Yau2}}
A Hom-Lie triple system (Hom-LTS) $(T, [\cdot,\cdot,\cdot], \alpha=(\alpha_1, \alpha_2))$ consists of an $\F$-vector space $T$, a trilinear map $[\cdot,\cdot,\cdot]: T\times T\times T\rightarrow T$, and linear maps $\alpha_i: T\rightarrow T$ for $i=1,2$, called twisted maps, such that for all $x,y,z,u,v\in T$,
\begin{enumerate}[(1)]
\item $[xxz]=0$,
\item $[xyz]+[yzx]+[zxy]=0$,
\item $[\alpha_1(u)\alpha_2(v)[xyz]] =[[uvx]\alpha_1(y)\alpha_2(z)] +[\alpha_1(x)[uvy]\alpha_2(z)]+[\alpha_1(x)\alpha_2(y)[uvz]]$.
\end{enumerate}
A Hom-LTS is said to be multiplicative if $\alpha_1=\alpha_2=\alpha$ and $\alpha([xyz])=[\alpha(x)\alpha(y)\alpha(z)]$, and denoted by $(T, [\cdot,\cdot,\cdot], \alpha)$.

A morphism $f: (T, [\cdot,\cdot,\cdot], \alpha=(\alpha_1, \alpha_2))\rightarrow (T', [\cdot,\cdot,\cdot]', \alpha'=(\alpha_1', \alpha_2'))$ of Hom-LTSs is a linear map satisfying $f([xyz])=[f(x)f(y)f(z)]'$ and $f\circ\alpha_i=\alpha_i'\circ f$ for $i=1,2$. An isomorphism is a bijective morphism.
\edefn

\bre
When the twisted maps $\alpha_i$ are both equal to the identity map, a Hom-LTS is a LTS. So LTSs are special examples of Hom-LTSs.
\ere

\beg{\rm\supercite{Yau2}}
\begin{enumerate}[(i)]
\item Let $V$ be a vector space over $\F$, $\langle , \rangle:V\times V\rightarrow \F$ be a symmetric bilinear form, and $\lambda\in\F$ be an arbitrary scalar. Suppose that $\alpha:V\rightarrow V$ is a linear map such that $\langle x,y\rangle=\langle \alpha(x),\alpha(y)\rangle$, for all $x, y\in V$. Define the triple product on $V$ by
$$[xyz]_\alpha=\lambda(\langle y,z\rangle\alpha(x)-\langle z,x\rangle\alpha(y)).$$
Then $(V, [\cdot,\cdot,\cdot]_\alpha, \alpha)$ is a multiplicative Hom-LTS.
\item Let $V=\F^{m\times n}$ be the vector space of $m\times n$ matrices with entries in $\F$ and $A^*$ denote the transpose of the matrix $A$. Suppose $\alpha: V\rightarrow V$ is a linear map such that $\alpha(AB)=\alpha(A)\alpha(B)$ and $\alpha(A^*)=\alpha(A)^*$. Then $(V, [\cdot,\cdot,\cdot]_\alpha, \alpha)$ is a multiplicative Hom-LTS with
    $$[ABC]_\alpha=\alpha(A)\alpha(B^*)\alpha(C)+ \alpha(C)\alpha(B^*)\alpha(A)- \alpha(B)\alpha(A^*)\alpha(C)-  \alpha(C)\alpha(A^*)\alpha(B).$$
\end{enumerate}\eeg

Following the representation theory introduced by Sheng for Hom-Lie algebras, and by Yamaguti for Lie triple systems, we generalize the notion of the representation to Hom-Lie triple systems.
\bdefn
Let $(T, [\cdot,\cdot,\cdot], \alpha)$ be a multiplicative Hom-LTS, $V$ an $\F$-vector space and $A\in End(V)$. $V$ is called a $(T, [\cdot,\cdot,\cdot], \alpha)$-module with respect to $A$ if there exists a bilinear map $\theta: T\times T\rightarrow End(V)$, $(a, b)\mapsto \theta(a, b)$ such that for all $a,b,c,d\in T$,
\begin{gather}
\theta(\alpha(a),\alpha(b))\circ A=A\circ \theta(a,b),\label{represen1}\\
\theta(\alpha(c),\alpha(d))\theta(a,b)\!-\! \theta(\alpha(b),\alpha(d))\theta(a,c)\!-\! \theta(\alpha(a),[bcd])\!\circ \!\! A\!\!+\! D(\alpha(b),\alpha(c))\theta(a,d)\!\!=\!\!0,\!\label{represen2}\\
\theta(\alpha(c),\alpha(d))D(a,b)\!-\! D(\alpha(a),\alpha(b))\theta(c,d)\!+\!\theta([abc],\alpha(d))\!\circ \!\!A\!+\!\theta(\alpha(c),[abd])\!\circ\!\! A\!=\!0,\label{represen3}
\end{gather}
where $D(a,b)=\theta(b,a)-\theta(a,b)$. Then $\theta$ is called the representation of $(T, [\cdot,\cdot,\cdot], \alpha)$ on $V$ with respect to $A$. In the case $\theta=0$, $V$ is called the trivial $(T, [\cdot,\cdot,\cdot], \alpha)$-module with respect to $A$.
\edefn
It can be concluded from (\ref{represen3}) that
\beq
D(\alpha(c),\alpha(d))D(a,b)\!-\!D(\alpha(a),\alpha(b))D(c,d)\! +\!D([abc],\alpha(d))\circ\! A\!+\!D(\alpha(c),[abd])\circ\! A\!=\!0.\label{represen4}
\eeq
In particular, let $V=T$, $A=\alpha$ and $\theta(x,y)(z)=[zxy]$. Then $D(x,y)(z)=[xyz]$ and (\ref{represen1}), (\ref{represen2}), (\ref{represen3}) hold. In this case $T$ is said to be the adjoint $(T, [\cdot,\cdot,\cdot], \alpha)$-module and $\theta$ is called the adjoint representation of $(T, [\cdot,\cdot,\cdot], \alpha)$ on itself with respect to $\alpha$.

As the case of general algebras, we give the semidirect product of a multiplicative Hom-LTS and its module.

\bprop
Let $\theta$ be a representation of a multiplicative Hom-LTS $(T, [\cdot,\cdot,\cdot], \alpha)$ on $V$ with respect to $A$. Then $T\oplus V$ has a structure as a multiplicative Hom-LTS.
\eprop
\bpf
We define the operation $[\cdot,\cdot,\cdot]_V: (T\oplus V)\times (T\oplus V)\times (T\oplus V)\rightarrow T\oplus V$ by
$$[(x,a),(y,b),(z,c)]_V=([xyz], \theta(y,z)(a)-\theta(x,z)(b)+D(x,y)(c));$$
and define the twisted map $\alpha+A: T\oplus V\rightarrow T\oplus V$ by
$$(\alpha+A)(x,a)=(\alpha(x), A(a)).$$
Using $[xxy]=0$, $[xyz]+[yzx]+[zxy]=0$ and $D(x,y)=\theta(y,x)-\theta(x,y)$, we obtain
$$[(x,a),(x,a),(y,b)]_V=(0,0),$$
$$[(x,a),(y,b),(z,c)]_V+[(y,b),(z,c)(x,a)]_V+[(z,c),(x,a),(y,b)]_V=(0, 0).$$
By (\ref{represen2}), (\ref{represen3}) and (\ref{represen4}), it follows that
\begin{align*}
&[[(x,a),(y,b),(u,c)]_V, (\alpha+A)(v,d), (\alpha+A)(w,e)]_V\\
&+ [(\alpha+A)(u,c), [(x,a),(y,b),(v,d)]_V, (\alpha+A)(w,e)]_V\\
&+ [(\alpha+A)(u,c), (\alpha+A)(v,d), [(x,a),(y,b),(w,e)]_V]_V\\
=&([[xyu]\alpha(v)\alpha(w)],\Omega_1) +([\alpha(u)[xyv]\alpha(w)],\Omega_2) +([\alpha(u)\alpha(v)[xyw]],\Omega_3)\\
=&([\alpha(x)\alpha(y)[uvw]],\Omega_4)\\
=&[(\alpha+A)(x,a),(\alpha+A)(y,b),[(u,c), (v,d), (w,e)]_V]_V,
\end{align*}
where
\begin{align*}
\Omega_1=&\theta(\alpha(v),\alpha(w))\theta(y,u)(a)-\theta(\alpha(v),\alpha(w))\theta(x,u)(b) +\theta(\alpha(v),\alpha(w))D(x,y)(c)\\
         &-\theta([xyu], \alpha(w))A(d)+D([xyu], \alpha(v))A(e),\\
\Omega_2=&-\theta(\alpha(u),\alpha(w))\theta(y,v)(a)+\theta(\alpha(u),\alpha(w))\theta(x,v)(b) +\theta([xyv], \alpha(w))A(c) \\
         &-\theta(\alpha(u),\alpha(w))D(x,y)(d)+D(\alpha(u), [xyv])A(e),\\
\Omega_3=&D(\alpha(u),\alpha(v))\theta(y,w)(a) -D(\alpha(u),\alpha(v))\theta(x,w)(b)+\theta(\alpha(v),[xyw])A(c)\\
         &-\theta(\alpha(u),[xyw])A(d)+D(\alpha(u),\alpha(v))D(x,y)(e),\\
\Omega_4=&\theta(\alpha(y),[uvw])A(a)-\theta(\alpha(x),[uvw])A(b) +D(\alpha(x),\alpha(y))\theta(v,w)(c)\\
         &-D(\alpha(x),\alpha(y))\theta(u,w)(d)+D(\alpha(x),\alpha(y))D(u,v)(e).
\end{align*}

Note that $\alpha+A$ is a linear map and using (\ref{represen1}), we have
\begin{align*}
&(\alpha+A)[(x,a),(y,b),(z,c)]_V\\
=&(\alpha+A)([xyz], \theta(y,z)(a)-\theta(x,z)(b)+D(x,y)(c))\\
=&(\alpha([xyz]), A\circ(\theta(y,z)(a)-\theta(x,z)(b)+D(x,y)(c)))\\
=&([\alpha(x)\alpha(y)\alpha(z)], \theta(\alpha(y),\alpha(z))A(a)-\theta(\alpha(x),\alpha(z))A(b)+D(\alpha(x),\alpha(y))A(c))\\
=&[(\alpha(x), A(a)),(\alpha(y), A(b)),(\alpha(z), A(c))]_V\\
=&[(\alpha+A)(x,a),(\alpha+A)(y,b),(\alpha+A)(z,c)]_V.
\end{align*}
Thus, $(T\oplus V, [\cdot,\cdot,\cdot]_V, \alpha+A)$ is a multiplicative Hom-LTS.
\epf
Let $\theta$ be a representation of $(T, [\cdot,\cdot,\cdot], \alpha)$ on $V$ with respect to $A$. If an $n$-linear map $f:\underbrace{T\times\cdots\times T}_{n~{\rm times}}\rightarrow V$ satisfies
\begin{gather}
A(f(x_1,\cdots,x_{n}))=f(\alpha(x_1),\cdots,\alpha(x_{n})),\label{cochain1}\\
f(x_1,\cdots,x,x,x_n)=0,\label{cochain2}\\
f(x_1,\cdots,x_{n-3},x,y,z) +f(x_1,\cdots,x_{n-3},y,z,x)+f(x_1,\cdots,x_{n-3},z,x,y)=0,\label{cochain3}
\end{gather}
then $f$ is called an $n$-Hom-cochain on $T$. Denote by $C_{\alpha,A}^n(T,V)$ the set of all $n$-Hom-cochains, $\forall n\geq1$.

\bdefn
For $n\geq1$, the coboundary operator $\delta^n_{hom}: C_{\alpha,A}^n(T,V) \rightarrow C_{\alpha,A}^{n+2}(T,V)$ is defined as follows.

For $f\in C_{\alpha,A}^{2n-1}(T,V), n=1,2,3,\cdots,$
\begin{align*}
     &\delta^{2n-1}_{hom}f(x_1,\cdots,x_{2n+1})\\
    =&\theta(\alpha^{n-1}(x_{2n}), \alpha^{n-1}(x_{2n+1}))f(x_1,\cdots,x_{2n-1})\\
     &-\theta(\alpha^{n-1}(x_{2n-1}), \alpha^{n-1}(x_{2n+1}))f(x_1,\cdots,x_{2n-2},x_{2n})\\
     &+\sum_{k=1}^n(-1)^{n+k}D(\alpha^{n-1}(x_{2k-1}), \alpha^{n-1}(x_{2k})) f(x_1,\cdots,\widehat{x_{2k-1}},\widehat{x_{2k}},\cdots, x_{2n+1})\\
     &+\sum_{k=1}^n\sum_{j=2k+1}^{2n+1}(-1)^{n+k+1} f(\alpha(x_1),\cdots,\widehat{x_{2k-1}},\widehat{x_{2k}},\cdots, [x_{2k-1}x_{2k}x_j],\cdots, \alpha(x_{2n+1}));
\end{align*}
for $f\in C_{\alpha,A}^{2n}(T,V), n=1,2,3,\cdots,$
\begin{align*}
     &\delta^{2n}_{hom}f(y,x_1,\cdots,x_{2n+1})\\
    =&\theta(\alpha^{n}(x_{2n}), \alpha^{n}(x_{2n+1}))f(y, x_1,\cdots,x_{2n-1})\\
     &-\theta(\alpha^{n}(x_{2n-1}), \alpha^{n}(x_{2n+1}))f(y, x_1,\cdots,x_{2n-2},x_{2n})\\
     &+\sum_{k=1}^n(-1)^{n+k}D(\alpha^{n}(x_{2k-1}), \alpha^{n}(x_{2k}))f(y, x_1,\cdots,\widehat{x_{2k-1}},\widehat{x_{2k}},\cdots, x_{2n+1})\\
     &+\sum_{k=1}^n\sum_{j=2k+1}^{2n+1}(-1)^{n+k+1} f(\alpha(y),\alpha(x_1),\cdots,\widehat{x_{2k-1}},\widehat{x_{2k}},\cdots, [x_{2k-1}x_{2k}x_j],\cdots, \alpha(x_{2n+1})),
\end{align*}
where the sign \textasciicircum ~ indicates that the element below must be omitted.
\edefn

It is not difficult to verify that $\delta^n_{hom}f$ satisfies (\ref{cochain2}) and (\ref{cochain3}). We claim that $\delta^n_{hom}f$ satisfies (\ref{cochain1}), so that the coboundary operator $\delta^n_{hom}$ is well-defined. In fact, it is sufficient to verify
$$A(\delta^{2n-1}_{hom}f(x_1,\cdots,x_{2n+1})) = \delta^{2n-1}_{hom}f(\alpha(x_1),\cdots,\alpha(x_{2n+1})).$$
However, by (\ref{represen1}) and $f$ being a $(2n-1)$-Hom-cochain, we have
\begin{align*}
 &\delta^{2n-1}_{hom}f(\alpha(x_1),\cdots,\alpha(x_{2n+1}))\\
=&\theta(\alpha^{n}(x_{2n}), \alpha^{n}(x_{2n+1}))f(\alpha(x_1),\cdots,\alpha(x_{2n-1}))\\
 &-\theta(\alpha^{n}(x_{2n-1}), \alpha^{n}(x_{2n+1}))f(\alpha(x_1),\cdots,\alpha(x_{2n-2}),\alpha(x_{2n}))\\
 &+\sum_{k=1}^n(-1)^{n+k}D(\alpha^{n}(x_{2k-1}), \alpha^{n}(x_{2k})) f(\alpha(x_1),\cdots,\widehat{\alpha(x_{2k-1})},\widehat{\alpha(x_{2k})}, \cdots, \alpha(x_{2n+1}))\\
 &+\sum_{k=1}^n\sum_{j=2k+1}^{2n+1}(-1)^{n+k+1} f(\alpha^2(x_1),\cdots,\widehat{x_{2k-1}},\widehat{x_{2k}},\cdots, \alpha([x_{2k-1}x_{2k}x_j]),\cdots, \alpha^2(x_{2n+1}))\\
=&\theta(\alpha^{n}(x_{2n}), \alpha^{n}(x_{2n+1}))A(f(x_1,\cdots,x_{2n-1}))\\
 &-\theta(\alpha^{n}(x_{2n-1}), \alpha^{n}(x_{2n+1}))A(f(x_1,\cdots,x_{2n-2},x_{2n}))\\
 &+\sum_{k=1}^n(-1)^{n+k}D(\alpha^{n}(x_{2k-1}), \alpha^{n}(x_{2k})) A(f(x_1,\cdots,\widehat{x_{2k-1}},\widehat{x_{2k}},\cdots, x_{2n+1}))\\
 &+\sum_{k=1}^n\sum_{j=2k+1}^{2n+1}(-1)^{n+k+1} A(f(\alpha(x_1),\cdots,\widehat{x_{2k-1}},\widehat{x_{2k}},\cdots, [x_{2k-1}x_{2k}x_j],\cdots, \alpha(x_{2n+1})))\\
=&A(\theta(\alpha^{n-1}(x_{2n}), \alpha^{n-1}(x_{2n+1}))f(x_1,\cdots,x_{2n-1}))\\
 &-A(\theta(\alpha^{n-1}(x_{2n-1}), \alpha^{n-1}(x_{2n+1}))f(x_1,\cdots,x_{2n-2},x_{2n}))\\
 &+A\left(\sum_{k=1}^n(-1)^{n+k}D(\alpha^{n-1}(x_{2k-1}), \alpha^{n-1}(x_{2k})) f(x_1,\cdots,\widehat{x_{2k-1}},\widehat{x_{2k}},\cdots, x_{2n+1})\right)\\
 &+A\left(\sum_{k=1}^n\sum_{j=2k+1}^{2n+1}(-1)^{n+k+1} f(\alpha(x_1),\cdots,\widehat{x_{2k-1}},\widehat{x_{2k}},\cdots, [x_{2k-1}x_{2k}x_j],\cdots, \alpha(x_{2n+1}))\right)\\
=&A(\delta^{2n-1}_{hom}f(x_1,\cdots,x_{2n+1})).
\end{align*}
\bthm
The coboundary operator $\delta^n_{hom}$ defined above satisfies $\delta^{n+2}_{hom}\delta^n_{hom}=0.$
\ethm
\bpf Observing the definition of the coboundary operator it follows immediately that $\delta^{2n+1}_{hom}\delta^{2n-1}_{hom}=0$ implies $\delta^{2n+2}_{hom}\delta^{2n}_{hom}=0$. Then we only need to prove $\delta^{2n+1}_{hom}\delta^{2n-1}_{hom}=0$ for $n=1,2,\cdots$.

Suppose that $f\in C_{\alpha,A}^{2n-1}(T,V), \forall n\geq 1$. Write $\delta^{2n-1}_{hom}$ and $\delta^{2n+1}_{hom}\delta^{2n-1}_{hom}$ as
$$\delta^{2n-1}_{hom}=\delta^{2n-1}_1+\delta^{2n-1}_2+\delta^{2n-1}_3 \text{\quad and \quad} \delta^{2n+1}_{hom}\delta^{2n-1}_{hom}=\sum_{i,j=1}^3\delta^{2n+1}_i\delta^{2n-1}_j,$$
where
\begin{align*}
\delta^{2n-1}_1f(x_1,\cdots,x_{2n+1}) = &\theta(\alpha^{n-1}(x_{2n}), \alpha^{n-1}(x_{2n+1}))f(x_1,\cdots,x_{2n-1})\\
&-\theta(\alpha^{n-1}(x_{2n-1}), \alpha^{n-1}(x_{2n+1}))f(x_1,\cdots,x_{2n-2},x_{2n}),
\end{align*}
\begin{align*}
&\delta^{2n-1}_2f(x_1,\cdots,x_{2n+1})\\ =& \sum_{k=1}^n(-1)^{n+k}
D(\alpha^{n-1}(x_{2k-1}), \alpha^{n-1}(x_{2k})) f(x_1,\cdots,\widehat{x_{2k-1}},\widehat{x_{2k}},\cdots, x_{2n+1}),
\end{align*}
\begin{align*}
&\delta^{2n-1}_3f(x_1,\cdots,x_{2n+1})\\ =& \sum_{k=1}^n\sum_{j=2k+1}^{2n+1}(-1)^{n+k+1}
f(\alpha(x_1),\cdots,\widehat{x_{2k-1}},\widehat{x_{2k}},\cdots, [x_{2k-1}x_{2k}x_j],\cdots, \alpha(x_{2n+1})).
\end{align*}

Let us compute first $\delta^{2n+1}_1\delta^{2n-1}_{hom}+ (\delta^{2n+1}_2+\delta^{2n+1}_3)\delta^{2n-1}_1$, and this is given by
\begin{align}
&(\delta^{2n+1}_1\delta^{2n-1}_{hom} +(\delta^{2n+1}_2+\delta^{2n+1}_3)\delta^{2n-1}_1) f(x_1,\cdots,x_{2n+3})\notag\\
=&\theta(\alpha^n(x_{2n+2}), \alpha^n(x_{2n+3}))\delta^{2n-1}_{hom}f(x_1,\cdots,x_{2n+1})\notag\\
 &-\theta(\alpha^n(x_{2n+1}), \alpha^n(x_{2n+3}))\delta^{2n-1}_{hom}f(x_1,\cdots,x_{2n}, x_{2n+2})\notag\\
 &+\sum_{k=1}^{n+1}(-1)^{n+k+1}D(\alpha^n(x_{2k-1}), \alpha^n(x_{2k}))\delta^{2n-1}_1 f(x_1,\cdots,\widehat{x_{2k-1}},\widehat{x_{2k}},\cdots, x_{2n+3})\notag\\
 &+\sum_{k=1}^{n+1}\sum_{j=2k+1}^{2n+3}(-1)^{n+k}\delta^{2n-1}_1 f(\alpha(x_1),\cdots,\widehat{x_{2k-1}},\widehat{x_{2k}},\cdots, [x_{2k-1}x_{2k}x_j],\cdots,  \alpha(x_{2n+3}))\notag\\
=&\theta(\alpha^n(x_{2n+2}), \alpha^n(x_{2n+3}))\notag\\
 &\relphantom{\theta}\cdot\Big(\theta(\alpha^{n-1}(x_{2n}), \alpha^{n-1}(x_{2n+1}))f(x_1,\cdots,x_{2n-1})\tag{a1}\\
 &\relphantom{\theta}-\theta(\alpha^{n-1}(x_{2n-1}), \alpha^{n-1}(x_{2n+1}))f(x_1,\cdots,x_{2n-2},x_{2n})\tag{a2}\\
 &\relphantom{\theta}+\sum_{k=1}^n(-1)^{n+k}D(\alpha^{n-1}(x_{2k-1}), \alpha^{n-1}(x_{2k})) f(x_1,\cdots,\widehat{x_{2k-1}},\widehat{x_{2k}},\cdots, x_{2n+1})\tag{b1}\\
 &\relphantom{\theta}+\!\sum_{k=1}^n\sum_{j=2k+1}^{2n+1}(-1)^{n+k+1} f(\alpha(x_1),\cdots,\widehat{x_{2k-1}},\widehat{x_{2k}},\cdots, [x_{2k-1}x_{2k}x_j],\cdots, \alpha(x_{2n+1}))\Big)\tag{c1}\\
 &-\theta(\alpha^n(x_{2n+1}), \alpha^n(x_{2n+3}))\notag\\
 &\relphantom{-\theta}\cdot\Big(\theta(\alpha^{n-1}(x_{2n}), \alpha^{n-1}(x_{2n+2}))f(x_1,\cdots,x_{2n-1})\tag{a3}\\
 &\relphantom{-\theta}-\theta(\alpha^{n-1}(x_{2n-1}), \alpha^{n-1}(x_{2n+2}))f(x_1,\cdots,x_{2n-2},x_{2n})\tag{a4}\\
 &\relphantom{-\theta}+\sum_{k=1}^n(-1)^{n+k}D(\alpha^{n-1}(x_{2k-1}), \alpha^{n-1}(x_{2k})) f(x_1,\cdots,\widehat{x_{2k-1}},\widehat{x_{2k}},\cdots, x_{2n}, x_{2n+2})\tag{b2}\\
 &\relphantom{-\theta}+\sum_{k=1}^n\sum_{j=2k+1}^{2n, 2n+2}\!(-1)^{n+k+1}\!f(\alpha(x_1),\!\cdots\!,\widehat{x_{2k-1}},\widehat{x_{2k}},\!\cdots\!, [x_{2k-1}x_{2k}x_j],\!\cdots\!, \alpha(x_{2n}), \alpha(x_{2n+2}))\!\Big)\tag{d1}\\
 &+\sum_{k=1}^{n}(-1)^{n+k+1}D(\alpha^n(x_{2k-1}), \alpha^n(x_{2k}))\notag\\
 &\relphantom{\sum}\cdot\Big(\theta(\alpha^{n-1}(x_{2n+2}), \alpha^{n-1}(x_{2n+3})) f(x_1,\cdots,\widehat{x_{2k-1}},\widehat{x_{2k}},\cdots, x_{2n+1})\tag{b3}\\
 &\relphantom{\sum}-\theta(\alpha^{n-1}(x_{2n+1}), \alpha^{n-1}(x_{2n+3})) f(x_1,\cdots,\widehat{x_{2k-1}},\widehat{x_{2k}},\cdots, x_{2n}, x_{2n+2})\Big)\tag{b4}\\
 &+D(\alpha^n(x_{2n+1}), \alpha^n(x_{2n+2}))\Big(\theta(\alpha^{n-1}(x_{2n}), \alpha^{n-1}(x_{2n+3}))f(x_1,\cdots, x_{2n-1})\tag{a5}\\
 &\relphantom{+D}-\theta(\alpha^{n-1}(x_{2n-1}), \alpha^{n-1}(x_{2n+3}))f(x_1,\cdots, x_{2n-2}, x_{2n})\Big)\tag{a6}\\
 &+\sum_{k=1}^{n}\sum_{j=2k+1}^{2n+1}(-1)^{n+k}\theta(\alpha^n(x_{2n+2}), \alpha^n(x_{2n+3}))f(\alpha(x_1),\cdots,\widehat{x_{2k-1}},\widehat{x_{2k}},\cdots, [x_{2k-1}x_{2k}x_j],\notag\\
 &\relphantom{+\sum}\cdots, \alpha(x_{2n+1}))\tag{c2}\\
 &+\sum_{k=1}^{n}(-1)^{n+k} \theta(\alpha^{n-1}[x_{2k-1}x_{2k}x_{2n+2}],\alpha^{n}(x_{2n+3})) f(\alpha(x_1),\cdots,\widehat{x_{2k-1}},\widehat{x_{2k}},\cdots, \alpha(x_{2n+1}))\tag{b5}\\
 &+\sum_{k=1}^{n}(-1)^{n+k} \theta(\alpha^{n}(x_{2n+2}),\alpha^{n-1}[x_{2k-1}x_{2k}x_{2n+3}]) f(\alpha(x_1),\cdots,\widehat{x_{2k-1}},\widehat{x_{2k}},\cdots, \alpha(x_{2n+1}))\tag{b6}\\
 &-\theta(\alpha^{n}(x_{2n}),\alpha^{n-1}[x_{2n+1}x_{2n+2}x_{2n+3}]) f(\alpha(x_1),\cdots,\alpha(x_{2n-1}))\tag{a7}\\
 &-\sum_{k=1}^{n}\sum_{j=2k+1}^{2n,2n+2}(-1)^{n+k}
 \theta(\alpha^n(x_{2n+1}), \alpha^n(x_{2n+3}))f(\alpha(x_1),\cdots,\widehat{x_{2k-1}},\widehat{x_{2k}},\cdots, [x_{2k-1}x_{2k}x_j],\notag\\
 &\relphantom{+\sum}\cdots, \alpha(x_{2n}), \alpha(x_{2n+2}))\tag{d2}\\
 &-\sum_{k=1}^{n}(\!-1)^{n+k}\theta (\!\alpha^{n-1}[x_{2k\!-\!1}x_{2k}x_{2n\!+\!1}],\alpha^{n}(x_{2n\!+\!3})) f(\!\alpha(x_1\!),\!\cdots\!,\widehat{x_{2k\!-\!1}},\widehat{x_{2k}},\!\cdots\!, \alpha(x_{2n}\!), \alpha(x_{2n\!+\!2}))\tag{b7}\\
 &-\sum_{k=1}^{n}(\!-1)^{n+k} \theta(\!\alpha^{n}(x_{2n\!+\!1}),\alpha^{n\!-\!1}[x_{2k\!-\!1}x_{2k}x_{2n\!+\!3}]) f(\!\alpha(x_1\!),\!\cdots\!,\widehat{x_{2k\!-\!1}},\widehat{x_{2k}},\!\cdots\!, \alpha(x_{2n}), \alpha(x_{2n\!+\!2}))\tag{b8}\\
 &+\theta(\alpha^{n}(x_{2n-1}),\alpha^{n-1}[x_{2n+1}x_{2n+2}x_{2n+3}]) f(\alpha(x_1),\cdots,\alpha(x_{2n-2}), \alpha(x_{2n})).\tag{a8}
\end{align}

By (\ref{represen2}) and (\ref{cochain1}), we have (a1)+$\cdots$+(a8)=0; by (\ref{represen3}) and (\ref{cochain1}), we conclude (b1)+$\cdots$+(b8)=0; it is obvious that (c1)+(c2)=(d1)+(d2)=0. Thus,
$$\delta^{2n+1}_1\delta^{2n-1}_{hom}+ (\delta^{2n+1}_2+\delta^{2n+1}_3)\delta^{2n-1}_1=0.$$

On the other hand,
\begin{align}
&(\delta^{2n+1}_2+\delta^{2n+1}_3)(\delta^{2n-1}_2+\delta^{2n-1}_3) f(x_1,\cdots,x_{2n+3})\notag\\
=&\sum_{k=1}^{n+1}(-1)^{n+k+1}D(\alpha^n(x_{2k-1}), \alpha^n(x_{2k})) (\delta^{2n-1}_2+\delta^{2n-1}_3) f(x_1,\cdots,\widehat{x_{2k-1}},\widehat{x_{2k}},\cdots, x_{2n+3})\notag\\
 &+\!\sum_{k=1}^{n+1}\sum_{j=2k+1}^{2n+3}(-1)^{n+k} (\delta^{2n-1}_2+\delta^{2n-1}_3) f(\alpha(x_1),\!\cdots\!,\!\widehat{x_{2k-1}},\widehat{x_{2k}},\!\cdots\!,\! [x_{2k-1}x_{2k}x_j],\!\cdots\!, \! \alpha(x_{2n+3}))\notag\\
=&-\sum_{1\leq i<k\leq n+1} (-1)^{k+i}D(\alpha^n(x_{2k-1}), \alpha^n(x_{2k})) D(\alpha^{n-1}(x_{2i-1}), \alpha^{n-1}(x_{2i}))f(x_1,\cdots,\widehat{x_{2i-1}},\widehat{x_{2i}},\notag\\ &\relphantom{+\sum}\cdots,\widehat{x_{2k-1}},\widehat{x_{2k}},\cdots, x_{2n+3})\tag{A1}\\
&+\sum_{1\leq k<i\leq n+1} (-1)^{k+i}D(\alpha^n(x_{2k-1}), \alpha^n(x_{2k})) D(\alpha^{n-1}(x_{2i-1}), \alpha^{n-1}(x_{2i}))f(x_1,\cdots,\widehat{x_{2k-1}},\widehat{x_{2k}},\notag\\ &\relphantom{+\sum}\cdots, \widehat{x_{2i-1}},\widehat{x_{2i}},\cdots, x_{2n+3})\tag{A2}\\
 &+\sum_{1\leq i<k\leq n+1}\sum_{j=2i+1\atop  j\neq 2k-1,2k}^{2n+3} (-1)^{k+i}D(\alpha^n(x_{2k-1}), \alpha^n(x_{2k}))f(\alpha(x_1),\cdots,\widehat{x_{2i-1}},\widehat{x_{2i}}, \cdots,\widehat{x_{2k-1}},\widehat{x_{2k}},\notag\\
 &\relphantom{+\sum} \cdots, [x_{2i-1}x_{2i}x_j], \cdots, \alpha(x_{2n+3}))\tag{B1}\\
 &-\sum_{1\leq k<i\leq n+1}\sum_{j=2i+1}^{2n+3} (-1)^{k+i}D(\alpha^n(x_{2k-1}), \alpha^n(x_{2k}))f(\alpha(x_1), \cdots,\widehat{x_{2k-1}},\widehat{x_{2k}},\cdots, \widehat{x_{2i-1}},\widehat{x_{2i}},\notag\\
 &\relphantom{+\sum}  \cdots, [x_{2i-1}x_{2i}x_j], \cdots, \alpha(x_{2n+3}))\tag{B2}\\
 &+\sum_{1\leq i<k\leq n+1}\sum_{j=2k+1}^{2n+3} (-1)^{k+i}D(\alpha^n(x_{2i-1}), \alpha^n(x_{2i}))f(\alpha(x_1),\cdots,\widehat{x_{2i-1}},\widehat{x_{2i}}, \cdots,\widehat{x_{2k-1}},\widehat{x_{2k}},\notag\\
 &\relphantom{+\sum} \cdots, [x_{2k-1}x_{2k}x_j], \cdots, \alpha(x_{2n+3}))\tag{B3}\\
 &-\sum_{1\leq k<i\leq n+1}\sum_{j=2k+1\atop  j\neq 2i-1,2i}^{2n+3} (-1)^{k+i}D(\alpha^n(x_{2i-1}), \alpha^n(x_{2i}))f(\alpha(x_1), \cdots,\widehat{x_{2k-1}},\widehat{x_{2k}},\cdots, \widehat{x_{2i-1}},\widehat{x_{2i}}, \notag\\
 &\relphantom{+\sum} \cdots, [x_{2k-1}x_{2k}x_j], \cdots, \alpha(x_{2n+3}))\tag{B4}\\
 &-\sum_{1\leq k<i\leq n+1} (-1)^{k+i}D(\alpha^{n-1}[x_{2k-1}x_{2k}x_{2i-1}], \alpha^n(x_{2i}))f(\alpha(x_1), \!\cdots\!,\widehat{x_{2k-1}},\widehat{x_{2k}},\!\cdots\!, \widehat{x_{2i-1}},\widehat{x_{2i}},\notag\\
 &\relphantom{+\sum}  \cdots, \alpha(x_{2n+3}))\tag{A3}\\
 &-\sum_{1\leq k<i\leq n+1} (-1)^{k+i}D(\alpha^n(x_{2i-1}),\alpha^{n-1}[x_{2k-1}x_{2k}x_{2i}])f(\alpha(x_1), \!\cdots\!,\widehat{x_{2k-1}},\widehat{x_{2k}},\!\cdots\!, \widehat{x_{2i-1}},\widehat{x_{2i}},\notag\\
 &\relphantom{+\sum}  \cdots, \alpha(x_{2n+3}))\tag{A4}\\
 &+\sum_{1\leq i<k\leq n+1}\sum_{2i<s<j\leq 2n+3\atop s\neq 2k-1,2k; 2k<j} (-1)^{k+i+1} f(\alpha^2(x_1), \!\cdots\!, \widehat{x_{2i-1}}, \widehat{x_{2i}}, \!\cdots\!, \widehat{x_{2k-1}},\widehat{x_{2k}}, \!\cdots\!, \alpha[x_{2i-1}x_{2i}x_s], \notag\\
 &\relphantom{+\sum} \cdots, \alpha[x_{2k-1}x_{2k}x_j], \cdots, \alpha^2(x_{2n+3}))\tag{C1}\\
 &+\sum_{1\leq i<k\leq n+1}\sum_{2k<j<s\leq 2n+3} (-1)^{k+i+1} f(\alpha^2(x_1), \!\cdots\!,\widehat{x_{2i-1}},\widehat{x_{2i}},\!\cdots\!, \widehat{x_{2k-1}},\widehat{x_{2k}}, \!\cdots\!, \alpha[x_{2k-1}x_{2k}x_j],\notag\\
 &\relphantom{+\sum} \cdots, \alpha[x_{2i-1}x_{2i}x_s], \cdots, \alpha^2(x_{2n+3}))\tag{C2}\\
 &-\sum_{1\leq k<i\leq n+1}\sum_{2i<s<j\leq 2n+3} (-1)^{k+i+1} f(\alpha^2(x_1), \!\cdots\!,\widehat{x_{2k-1}},\widehat{x_{2k}},\!\cdots\!, \widehat{x_{2i-1}},\widehat{x_{2i}}, \!\cdots\!, \alpha[x_{2i-1}x_{2i}x_s],\notag\\
 &\relphantom{+\sum} \cdots, \alpha[x_{2k-1}x_{2k}x_j], \cdots, \alpha^2(x_{2n+3}))\tag{C3}\\
 &-\sum_{1\leq k<i\leq n+1}\sum_{2k<j<s\leq 2n+3\atop j\neq 2i-1,2i; 2i<s} (-1)^{k+i+1} f(\alpha^2(x_1), \!\cdots\!, \widehat{x_{2k-1}}, \widehat{x_{2k}}, \!\cdots\!, \widehat{x_{2i-1}}, \widehat{x_{2i}}, \!\cdots\!, \alpha[x_{2k-1}x_{2k}x_j], \notag\\
 &\relphantom{+\sum} \cdots, \alpha[x_{2i-1}x_{2i}x_s], \cdots, \alpha^2(x_{2n+3}))\tag{C4}\\
 &+\sum_{1\leq i<k\leq n+1}\sum_{j=2k+1}^{2n+3} (-1)^{k+i+1} f(\alpha^2(x_1), \cdots,\widehat{x_{2i-1}},\widehat{x_{2i}},\cdots, \widehat{x_{2k-1}},\widehat{x_{2k}}, \notag\\
 &\relphantom{+\sum} \cdots, [\alpha(x_{2i-1})\alpha(x_{2i})[x_{2k-1}x_{2k}x_j]],\cdots, \alpha^2(x_{2n+3}))\tag{D1}\\
 &-\sum_{1\leq k<i\leq n+1}\sum_{j=2i+1}^{2n+3} (-1)^{k+i+1} f(\alpha^2(x_1), \cdots,\widehat{x_{2k-1}},\widehat{x_{2k}},\cdots, \widehat{x_{2i-1}},\widehat{x_{2i}}, \notag\\
 &\relphantom{+\sum} \cdots, [\alpha(x_{2i-1})\alpha(x_{2i})[x_{2k-1}x_{2k}x_j]],\cdots, \alpha^2(x_{2n+3}))\tag{D2}\\
 &-\sum_{1\leq k<i\leq n+1}\sum_{s=2i+1}^{2n+3} (-1)^{k+i+1} f(\alpha^2(x_1), \cdots,\widehat{x_{2k-1}},\widehat{x_{2k}},\cdots, \widehat{x_{2i-1}},\widehat{x_{2i}},\notag\\
 &\relphantom{+\sum}  \cdots, [[x_{2k-1}x_{2k}x_{2i-1}]\alpha(x_{2i})\alpha(x_{s})],\cdots, \alpha^2(x_{2n+3}))\tag{D3}\\
 &-\sum_{1\leq k<i\leq n+1}\sum_{s=2i+1}^{2n+3} (-1)^{k+i+1} f(\alpha^2(x_1), \cdots,\widehat{x_{2k-1}},\widehat{x_{2k}},\cdots, \widehat{x_{2i-1}},\widehat{x_{2i}}, \notag\\
 &\relphantom{+\sum} \cdots, [\alpha(x_{2i-1})[x_{2k-1}x_{2k}x_{2i}]\alpha(x_{s})],\cdots, \alpha^2(x_{2n+3}))\tag{D4}\\
=&0,\notag
\end{align}
since it is straightforward to verify that the sum of terms labeled with the same letter vanishes(e.g. (A1) + $\cdots$ + (A4)=0). Therefore, $\delta^{2n+1}_{hom}\delta^{2n-1}_{hom}=0$, and the proof is completed.
\epf

The map $f\in C_{\alpha, A}^n(T,V)$ is called an $n$-Hom-cocycle if $\delta^{n}_{hom}f=0$. We denote by $Z_{\alpha, A}^n(T,V)$ the subspace spanned by $n$-Hom-cocycles and $B_{\alpha, A}^n(T,V)=\delta^{n-2}_{hom}C_{\alpha, A}^{n-2}(T,V)$. Since $\delta^{n+2}_{hom}\delta^n_{hom}=0$, $B_{\alpha, A}^n(T,V)$ is a subspace of $Z_{\alpha, A}^n(T,V)$. Hence we can define a cohomology space $H_{\alpha, A}^n(T,V)$ of $(T, [\cdot,\cdot,\cdot], \alpha)$ as the factor space $Z_{\alpha, A}^n(T,V)/B_{\alpha, A}^n(T,V)$.

\bre
When $(T, [\cdot,\cdot,\cdot], \alpha)$ is a LTS, that is, $\alpha=\id_T$, the cohomology theory for Hom-LTSs above is actually the cohomology theory for LTSs in \cite{Yamaguti}.
\ere

\section{Central extensions of Hom-Lie triple systems}

Suppose that $(T, [\cdot,\cdot,\cdot], \alpha)$ is a multiplicative Hom-LTS, and $V$ is a trivial $(T, [\cdot,\cdot,\cdot], \alpha)$-module with respect to $\alpha_V$. Obviously, $(V, 0, \alpha_V)$ is an abelian multiplicative Hom-LTS with the trivial product. A multiplicative Hom-LTS $(T_C, [\cdot,\cdot,\cdot]_C, \alpha_C)$ is called \textbf{a central extension} of $(T, [\cdot,\cdot,\cdot], \alpha)$ by $(V, 0, \alpha_V)$ if there is a diagram with exact rows of Hom-LTSs
$$\xymatrix
{0\ar[r]& V\ar[r]^\iota\ar[d]^{\alpha_V}& T_C\ar[r]^\pi\ar[d]^{\alpha_C}& T\ar[r]\ar[d]^{\alpha}\ar@/^/[l]^{s}& 0\\
0\ar[r]& V\ar[r]^\iota& T_C\ar[r]^\pi& T\ar[r]\ar@/^/[l]^{s}& 0}$$
such that $\alpha_C\circ \iota=\iota\circ \alpha_V$ and $\alpha\circ\pi=\pi\circ\alpha_C$, where $s$ is a linear map satisfying $\pi s=\id_T$ and $\alpha_C\circ s=s\circ \alpha$, and $\iota(V)$ is contained in $Z(T_C)=\{x\in T_C \ |\ [x,T_C,T_C]_C=0\}$, the center of $T_C$.
If $(T_C', [\cdot,\cdot,\cdot]_C', \alpha_C')$ is another central extension  of $(T, [\cdot,\cdot,\cdot], \alpha)$ by $(V, 0, \alpha_V)$ such that there is an isomorphism $\varphi:(T, [\cdot,\cdot,\cdot], \alpha)\rightarrow (T_C', [\cdot,\cdot,\cdot]_C', \alpha_C')$, and the diagram
$$\xymatrix
{0\ar[r]& V\ar[r]^\iota\ar[d]^{\id_V}& T_C\ar[r]^\pi\ar[d]^{\varphi}& T\ar[r]\ar[d]^{\id_T}& 0\\
0\ar[r]& V\ar[r]^{\iota'}& T_C'\ar[r]^{\pi'}& T\ar[r]& 0}$$
commutes, then the two central extensions $(T_C, [\cdot,\cdot,\cdot]_C, \alpha_C)$ and $(T_C', [\cdot,\cdot,\cdot]_C', \alpha_C')$ are said to be \textbf{equivalent}.
\bthm
There is a one-to-one correspondence between equivalent classes of central extensions of $(T, [\cdot,\cdot,\cdot], \alpha)$ by $(V, 0, \alpha_V)$ and $H_{\alpha, \alpha_V}^3(T,V)$.
\ethm
\bpf
First we will construct a one-to-one correspondence between central extensions of $(T, [\cdot,\cdot,\cdot], \alpha)$ by $(V, 0, \alpha_V)$ and $Z_{\alpha, \alpha_V}^3(T,V)$.

Suppose that $(T_C, [\cdot,\cdot,\cdot]_C, \alpha_C)$ is a central extension of $(T, [\cdot,\cdot,\cdot], \alpha)$ by $(V, 0, \alpha_V)$. Then we have the following diagram
$$\xymatrix
{0\ar[r]& V\ar[r]^\iota\ar[d]^{\alpha_V}& T_C\ar[r]^\pi\ar[d]^{\alpha_C}& T\ar[r]\ar[d]^{\alpha}\ar@/^/[l]^{s}& 0\\
0\ar[r]& V\ar[r]^\iota& T_C\ar[r]^\pi& T\ar[r]\ar@/^/[l]^{s}& 0}$$
with $\alpha_C\circ \iota=\iota\circ \alpha_V$, $\alpha\circ\pi=\pi\circ\alpha_C$, and $s$ being a linear map satisfying $\pi s=\id_T$ and $\alpha_C\circ s=s\circ \alpha$. We will find an element in $Z_{\alpha, \alpha_V}^3(T,V)$.

For $x,y,z\in T$, since $\pi[s(x)s(y)s(z)]_C-\pi s[xyz]=[\pi s(x)\pi s(y)\pi s(z)]-[xyz]=0$, it follows that $[s(x)s(y)s(z)]_C-s[xyz]\in \Ker \pi=\iota(V)$. Define a trilinear map $g:T\times T\times T\rightarrow V$ by
$$\iota g(x,y,z)=[s(x)s(y)s(z)]_C-s[xyz].$$
Since $\iota$ is injective, $g$ is well-defined, and it follows from $\iota(V)\subseteq Z(T_C)$ that
$$[[s(x)s(y)s(z)]_C,u,v]_C=[s[xyz],u,v]_C, \forall u,v\in T_C.$$
Note that $g$ satisfies $g(x,x,y)=0, g(x,y,z)+g(y,z,x)+g(z,x,y)=0$ and
\begin{align*}
\iota g(\alpha(x),\alpha(y)\alpha(z))=&[s\alpha(x),s\alpha(y),s\alpha(z)]_C-s[\alpha(x)\alpha(y)\alpha(z)]\\
=&[\alpha_C s(x),\alpha_C s(y),\alpha_C s(z)]_C-\alpha_C s[xyz]\\
=&\alpha_C([s(x)s(y)s(z)]_C-s[xyz])\\
=&\alpha_C\iota g(x,y,z)=\iota\alpha_V g(x,y,z).
\end{align*}
Hence $g\in C_{\alpha, \alpha_V}^3(T,V)$. Moreover, $g\in Z_{\alpha, \alpha_V}^3(T,V)$ since
\begin{align*}
&\iota(\delta_{hom}^3 g)(x_1,x_2,x_3,x_4,x_5)\\
=&\iota(g([x_1x_2x_3],\alpha(x_4),\alpha(x_5))+g(\alpha(x_3),[x_1x_2x_4],\alpha(x_5))+g(\alpha(x_3),\alpha(x_4),[x_1x_2x_5])\\
 &-g(\alpha(x_1),\alpha(x_2),[x_3x_4x_5]))\\
=&[s[x_1x_2x_3],s\alpha(x_4),s\alpha(x_5)]_C-s[[x_1x_2x_3]\alpha(x_4)\alpha(x_5)]\\
 &+[s\alpha(x_3),s[x_1x_2x_4],s\alpha(x_5)]_C-s[\alpha(x_3)[x_1x_2x_4]\alpha(x_5)]\\
 &+[s\alpha(x_3),s\alpha(x_4),s[x_1x_2x_5]]_C-s[\alpha(x_3)\alpha(x_4)[x_1x_2x_5]]\\
 &-[s\alpha(x_1),s\alpha(x_2),s[x_3x_4x_5]]_C+s[\alpha(x_1)\alpha(x_2)[x_3x_4x_5]]\\
=&[[s(x_1)s(x_2)s(x_3)]_C,\alpha_C s(x_4),\alpha_C s(x_5)]_C+[\alpha_C s(x_3),[s(x_1)s(x_2)s(x_4)]_C,\alpha_C s(x_5)]_C\\
 &+[\alpha_C s(x_3),\alpha_C s(x_4),[s(x_1)s(x_2)s(x_5)]_C]_C-[\alpha_C s(x_1),\alpha_C s(x_2),[s(x_3)s(x_4)s(x_5)]_C]_C\\
=&0.
\end{align*}

Conversely, let $g\in Z_{\alpha, \alpha_V}^3(T,V)$ and $T_C=T\oplus V$ with
$$[(x,a),(y,b),(z,c)]_C=([xyz],g(x,y,z));\quad \alpha_C(x,a)=(\alpha(x),\alpha_V(a)).$$
Then $\alpha_C$ is linear and
\begin{align*}
&\alpha_C[(x,a),(y,b),(z,c)]_C=\alpha_C([xyz],g(x,y,z))=(\alpha[xyz],\alpha_Vg(x,y,z))\\
=&([\alpha(x)\alpha(y)\alpha(z)], g(\alpha(x),\alpha(y),\alpha(z)))=[(\alpha(x),\alpha_V(a)),(\alpha(y),\alpha_V(b)),(\alpha(z),\alpha_V(c))]_C\\
=&[\alpha_C(x,a),\alpha_C(y,b),\alpha_C(z,c)]_C.
\end{align*}
Note that
\begin{align*}
&[\alpha_C(x,a),\alpha_C(y,b),[(u,c),(v,d),(w,e)]_C]_C\\
=&[(\alpha(x),\alpha_V(a)),(\alpha(y),\alpha_V(b)),([uvw],g(u,v,w))]_C\\
=&([\alpha(x)\alpha(y)[uvw]], g(\alpha(x),\alpha(y),[uvw]))\\
=&([[xyu]\alpha(v)\alpha(w)],g([xyu],\alpha(v),\alpha(w)))\\
 &+([\alpha(u)[xyv]\alpha(w)],g(\alpha(u),[xyv],\alpha(w)))\\
 &+([\alpha(u)\alpha(v)[xyw]],g(\alpha(u),\alpha(v),[xyw]))\\
=&[[(x,a),(y,b),(u,c)]_C,\alpha_C(v,d),\alpha_C(w,e)]_C+[\alpha_C(u,c),[(x,a),(y,b),(v,d)]_C,\alpha_C(w,e)]_C\\
 &+[\alpha_C(u,c),\alpha_C(v,d),[(x,a),(y,b),(w,e)]_C]_C.
\end{align*}
Thus, $(T_C, [\cdot,\cdot,\cdot]_C, \alpha_C)$ is a multiplicative Hom-LTS.

Define $\iota:V\rightarrow T_C$ by $\iota(a)=(0,a)$, $\pi:T_C\rightarrow T$ by $\pi(x,a)=x$ and $s: T\rightarrow T_C$ by $s(x)=(x,0)$. Then
$$\alpha_C\circ \iota(a)=\alpha_C(0,a)=(0,\alpha_V(a))=\iota\circ\alpha_V(a),$$
$$\pi\circ\alpha_C(x,a)=\pi(\alpha(x),\alpha_V(a))=\alpha(x)=\alpha\circ\pi(x,a),$$
$$\pi s=\id_T,\quad \alpha_C s(x)=\alpha_C(x,0)=(\alpha(x),0)=s\alpha(x).$$
It is clear $\iota(V)\subseteq Z(T_C)$. Therefore, $(T_C, [\cdot,\cdot,\cdot]_C, \alpha_C)$ is a central extension of $(T, [\cdot,\cdot,\cdot], \alpha)$ by $(V, 0, \alpha_V)$.

Suppose that $(T_C, [\cdot,\cdot,\cdot]_C, \alpha_C)$ and $(T_C', [\cdot,\cdot,\cdot]_C', \alpha_C')$ are equivalent central extensions of $(T, [\cdot,\cdot,\cdot], \alpha)$ by $(V, 0, \alpha_V)$. Then we have the following diagram
$$\xymatrix
{0\ar[r]& V\ar[r]^\iota\ar[d]^{\id_V}& T_C\ar[r]^\pi\ar[d]^{\varphi}& T\ar[r]\ar[d]^{\id_T}\ar@/^/[l]^{s}& 0\\
0\ar[r]& V\ar[r]^{\iota'}& T_C'\ar[r]^{\pi'}& T\ar[r]\ar@/^/[l]^{s'}& 0}$$
such that $\varphi\circ\iota=\iota'$ and $\pi=\pi'\circ\varphi$, with $\varphi$ being an isomorphism and $\pi s=\pi's'=\id_T$. Let $g,g'$ be their corresponding 3-Hom-cocycles constructed as above, respectively. Then
$$\iota g(x,y,z)=[s(x)s(y)s(z)]_C-s[xyz],$$
$$\iota' g'(x,y,z)=[s'(x)s'(y)s'(z)]_C'-s'[xyz],$$
$$\iota'g(x,y,z)=\varphi\iota g(x,y,z)=\varphi[s(x)s(y)s(z)]_C-\varphi s[xyz].$$
We claim that $g-g'\in B_{\alpha, \alpha_V}^3(T,V)$. In fact, using
$$\pi's'(x)-\pi'\varphi s(x)=x-\pi s(x)=0,$$
we define a linear map $f:T\rightarrow V$ by $\iota'f(x)=s'(x)-\varphi s(x)$, for all $x\in T$. Then
\begin{align*}
\iota'f\alpha(x)&=s' \alpha(x)-\varphi s \alpha(x)=\alpha_C' s'(x)-\varphi\alpha_C s(x)\\
                &=\alpha_C' s'(x)-\alpha_C'\varphi s(x)=\alpha_C'\iota'f(x)=\iota'\alpha_V f(x),
\end{align*}
which implies $f\in C_{\alpha, \alpha_V}^1(T,V)$. Since $s'(x)-\varphi s(x)=\iota'f(x)\in Z(T_C')$,
$$[s'(x)s'(y)s'(z)]_C'=[\varphi s(x)\varphi s(y)\varphi s(z)]_C'=\varphi[s(x)s(y)s(z)]_C.$$
Then
$$\iota'(g'-g)(x,y,z)=-\iota'f([xyz])=\iota'(\delta_{hom}^1 f)(x,y,z),$$
so $g'-g=\delta_{hom}^1 f\in B_{\alpha, \alpha_V}^3(T,V)$.

Suppose $g,g'\in Z_{\alpha, \alpha_V}^3(T,V)$ and $g'-g\in B_{\alpha, \alpha_V}^3(T,V)$, i.e., there is $f\in C_{\alpha, \alpha_V}^1(T,V)$ satisfying $g'-g=\delta_{hom}^1 f$. Then $(g'-g)(x,y,z)=-f([xyz])$. Let $(T_C, [\cdot,\cdot,\cdot]_C, \alpha_C)$ and $(T_C', [\cdot,\cdot,\cdot]_C', \alpha_C)$ be two central extensions of $(T, [\cdot,\cdot,\cdot], \alpha)$ by $(V, 0, \alpha_V)$ defined as above with respect to $g$ and $g'$, respectively. Then $\iota(a)=\iota'(a)=(0,a)$ and $\pi(x,a)=\pi'(x,a)=x$. Consider a linear map
\begin{eqnarray*}
\varphi: (T_C, [\cdot,\cdot,\cdot]_C, \alpha_C) &\longrightarrow& (T_C', [\cdot,\cdot,\cdot]_C', \alpha_C)\\
         (x,a)                                  &\longmapsto    & (x,a-f(x)).
\end{eqnarray*}
Then $\varphi\iota(a)=\iota'(a)$ and $\pi'\varphi(x,a)=\pi'(x,a-f(x))=x=\pi(x,a).$ Hence we obtain the following commutative diagram
$$\xymatrix
{0\ar[r]& V\ar[r]^\iota\ar[d]^{\id_V}& T_C\ar[r]^\pi\ar[d]^{\varphi}& T\ar[r]\ar[d]^{\id_T}& 0\\
0\ar[r]& V\ar[r]^{\iota'}& T_C\ar[r]^{\pi'}& T\ar[r]& 0.}$$
Now it is sufficient to prove that $\varphi$ is an isomorphism.

If $\varphi(x,a)=\varphi(\tilde{x},\tilde{a})$, it follows that $(x, a-f(x))=(\tilde{x}, \tilde{a}-f(\tilde{x}))$, that is, $x=\tilde{x}$ and $a-f(x)=\tilde{a}-f(\tilde{x})$, then $a=\tilde{a}$, and so $\varphi$ is injective; $\varphi$ is obviously surjective. Note that
\begin{align*}
\varphi\alpha_C(x,a)
&=\varphi(\alpha(x),\alpha_V(a))=(\alpha(x),\alpha_V(a)-f\alpha(x)) =(\alpha(x),\alpha_V(a)-\alpha_Vf(x))\\
&=\alpha_C'(x, a-f(x))=\alpha_C'\varphi(x,a)
\end{align*}
and
\begin{align*}
\varphi[(x,a),(y,b),(z,c)]_C&=\varphi([xyz],g(x,y,z))=([xyz],g(x,y,z)-f([xyz]))\\
                          &=([xyz],g'(x,y,z))=[(x,a-f(x)),(y,b-f(y)),(z,c-f(z))]_C'\\
                          &=[\varphi(x,a),\varphi(y,b),\varphi(z,c)]_C'.
\end{align*}

Therefore, $(T_C, [\cdot,\cdot,\cdot]_C, \alpha_C)$ and $(T_C', [\cdot,\cdot,\cdot]_C', \alpha_C')$ are equivalent central extensions of $(T, [\cdot,\cdot,\cdot], \alpha)$ by $(V, 0, \alpha_V)$.
\epf

\section{1-parameter formal deformations of multiplicative Hom-Lie triple systems}

Let $(T, [\cdot,\cdot,\cdot], \alpha)$ be a multiplicative Hom-LTS and $\F[[t]]$ be the ring of formal power series over $\F$. Suppose that $T[[t]]$ is the set of formal power series over $T$. Then for an $\F$-trilinear map $f:T\times T\times T\rightarrow T$, it is natural to extend it to be an $\F[[t]]$-trilinear map $f:T[[t]]\times T[[t]]\times T[[t]]\rightarrow T[[t]]$ by
$$f\left(\sum_{i\geq0}x_it^i,\sum_{j\geq0}y_jt^j,\sum_{k\geq0}z_kt^k\right)=\sum_{i,j,k\geq0}f(x_i,y_j,z_k)t^{i+j+k}.$$
\bdefn
Let $(T, [\cdot,\cdot,\cdot], \alpha)$ be a multiplicative Hom-LTS over $\F$. A 1-parameter formal deformation of $(T, [\cdot,\cdot,\cdot], \alpha)$ is a formal power series $d_t:T[[t]]\times T[[t]]\times T[[t]]\rightarrow T[[t]]$ of the form
$$d_t(x,y,z)=\sum_{i\geq0}d_i(x,y,z)t^i=d_0(x,y,z)+d_1(x,y,z)t+d_2(x,y,z)t^2+\cdots,$$
where each $d_i$ is an $\F$-trilinear map $d_i:T\times T\times T\rightarrow T$ (extended to be $\F[[t]]$-trilinear) and $d_0(x,y,z)=[xyz]$, such that the following identities hold
\begin{gather}
d_t(\alpha(x),\alpha(y),\alpha(z))=\alpha\circ d_t(x,y,z),\label{deq1}\\
d_t(x,x,y)=0,\label{deq2}\\
d_t(x,y,z)+d_t(y,z,x)+d_t(z,x,y)=0,\label{deq3}\\
\begin{aligned}
d_t(\alpha(u),\alpha(v),d_t(x,y,z))=&d_t(d_t(u,v,x),\alpha(y),\alpha(z))\\
                                    &\!+\!d_t(\alpha(x),d_t(u,v,y),\alpha(z))\!+\!d_t(\alpha(x),\alpha(y),d_t(u,v,z)).\label{deq4}
\end{aligned}
\end{gather}
Conditions (\ref{deq1})-(\ref{deq4}) are called the deformation equations of a multiplicative Hom-LTS.
\edefn

Note that $T[[t]]$ is a module over $\F[[t]]$ and $d_t$ defines the trilinear multiplication on $T[[t]]$ such that $T_t=(T[[t]],d_t,\alpha)$ is a multiplicative Hom-LTS. Now we investigate the deformation equations (\ref{deq1})-(\ref{deq4}).

Conditions (\ref{deq1})-(\ref{deq3}) are equivalent to the following equations
\begin{gather}
d_i(\alpha(x),\alpha(y),\alpha(z))=\alpha\circ d_i(x,y,z),\label{deq5}\\
d_i(x,x,y)=0,\label{deq6}\\
d_i(x,y,z)+d_i(y,z,x)+d_i(z,x,y)=0,\label{deq7}
\end{gather}
respectively, for $i=0,1,2,\cdots.$ The condition (\ref{deq4}) can be expressed as
\begin{align*}
&\sum_{i,j\geq0}d_i(\alpha(u),\alpha(v),d_j(x,y,z))\\
=&\sum_{i,j\geq0}\!d_i(d_j(u,v,x),\alpha(y),\alpha(z)) \!+\!\!\!\sum_{i,j\geq0}\!d_i(\alpha(x),d_j(u,v,y)\alpha(z)) \!+\!\!\!\sum_{i,j\geq0}\!d_i(\alpha(x),\alpha(y),d_j(u,v,z)).
\end{align*}
Then
\begin{align*}
\sum_{i+j=n}\Big(&d_i(d_j(u,v,x),\alpha(y),\alpha(z))+d_i(\alpha(x),d_j(u,v,y)\alpha(z))\\
&+d_i(\alpha(x),\alpha(y),d_j(u,v,z))-d_i(\alpha(u),\alpha(v),d_j(x,y,z))\Big)=0,\quad \forall n=0,1,2\cdots.
\end{align*}

For two $\F$-trilinear maps $f,g:T\times T\times T\rightarrow T$ (extended to be $\F[[t]]$-trilinear), define a map $f\circ_{\alpha}g:T[[t]]\times T[[t]]\times T[[t]]\times T[[t]]\times T[[t]]\rightarrow T[[t]]$ by
\begin{align*}
f\circ_{\alpha}g(u,v,x,y,z)=&f(g(u,v,x),\alpha(y),\alpha(z))+f(\alpha(x),g(u,v,y),\alpha(z))\\
&+f(\alpha(x),\alpha(y),g(u,v,z))-f(\alpha(u),\alpha(v),g(x,y,z)).
\end{align*}
Then the deformation equation (\ref{deq4}) can be written as
$$\sum_{i+j=n}d_i\circ_{\alpha}d_j=0.$$

For $n=1$, $d_0\circ_{\alpha}d_1+d_1\circ_{\alpha}d_0=0.$

For $n\geq2$, $-(d_0\circ_{\alpha}d_n+d_n\circ_{\alpha}d_0)=d_1\circ_{\alpha}d_{n-1}+d_2\circ_{\alpha}d_{n-2}+\cdots+d_{n-1}\circ_{\alpha}d_1.$

Section 2 says that $T$ is the adjoint $(T, [\cdot,\cdot,\cdot], \alpha)$-module by setting $\theta(x,y)(z)=[zxy]$ and $A=\alpha$. In this case, by (\ref{deq5})-(\ref{deq7}) it follows that $d_i\in C_{\alpha,\alpha}^3(T,T)$. It can also be verified that $d_i\circ_{\alpha}d_j\in C_{\alpha,\alpha}^5(T,T)$. In general, if $f, g\in C_{\alpha,\alpha}^3(T,T)$, then $f\circ_{\alpha}g\in C_{\alpha,\alpha}^5(T,T)$. Note that the coboundary operator can be written for lower orders as
\begin{align*}
\delta_{hom}^1f(x_1,x_2,x_3)=&[f(x_1)x_2x_3]+[x_1f(x_2)x_3]+[x_1x_2f(x_3)]-f([x_1x_2x_3]),\\
\delta_{hom}^3f(x_1,x_2,x_3,x_4,x_5)=&[f(x_1,x_2,x_3)\alpha(x_4)\alpha(x_5)]+[\alpha(x_3)f(x_1,x_2,x_4)\alpha(x_5)]\\
                                     &+[\alpha(x_3)\alpha(x_4)f(x_1,x_2,x_5)]-[\alpha(x_1)\alpha(x_2)f(x_3,x_4,x_5)]\\
                                     &+f([x_1x_2x_3],\alpha(x_4),\alpha(x_5))+f(\alpha(x_3),[x_1x_2x_4],\alpha(x_5))\\
                                     &+f(\alpha(x_3),\alpha(x_4),[x_1x_2x_5])-f(\alpha(x_1),\alpha(x_2),[x_3x_4x_5]),
\end{align*}
which implies $\delta_{hom}^3d_n=d_0\circ_{\alpha}d_n+d_n\circ_{\alpha}d_0,$ for $n=0,1,2\cdots.$ Hence the deformation equation (\ref{deq4}) can be rewritten as
\begin{gather*}
\delta_{hom}^3d_1=0,\\
-\delta_{hom}^3d_n=d_1\circ_{\alpha}d_{n-1}+d_2\circ_{\alpha}d_{n-2}+\cdots+d_{n-1}\circ_{\alpha}d_1.
\end{gather*}
Then $d_1$ is a 3-Hom-cocycle and called the \textbf{infinitesimal} of $d_t$.
\bdefn
Let $(T, [\cdot,\cdot,\cdot], \alpha)$ be a multiplicative Hom-LTS. Suppose that $d_t(x,y,z)=\sum_{i\geq0}d_i(x,y,z)t^i$ and $d_t'(x,y,z)=\sum_{i\geq0}d_i'(x,y,z)t^i$ are two 1-parameter formal deformations of $(T, [\cdot,\cdot,\cdot], \alpha)$. They are called equivalent, denoted by $d_t\sim d_t'$, if there is a formal isomorphism of $\F[[t]]$-modules
$$\phi_t(x,y,z)=\sum_{i\geq0}\phi_i(x,y,z)t^i:(T[[t]],d_t,\alpha)\longrightarrow (T[[t]],d_t',\alpha),$$
where each $\phi_i:T\rightarrow T$ is an $\F$-linear map (extended to be $\F[[t]]$-linear) and $\phi_0=\id_T$, satisfying
\begin{gather*}
\phi_t\circ\alpha=\alpha\circ\phi_t,\\
\phi_t\circ d_t(x,y,z)=d_t'(\phi_t(x),\phi_t(y),\phi_t(z)).
\end{gather*}
When $d_1=d_2=\cdots=0$, $d_t=d_0$ is said to be the null deformation. A 1-parameter formal deformation $d_t$ is called trivial if $d_t\sim d_0$.
\edefn
\bthm
Let $d_t(x,y,z)=\sum_{i\geq0}d_i(x,y,z)t^i$ and $d_t'(x,y,z)=\sum_{i\geq0}d_i'(x,y,z)t^i$ be equivalent 1-parameter formal deformations of $(T, [\cdot,\cdot,\cdot], \alpha)$. Then $d_1$ and $d_1'$ belong to the same cohomology class in $H_{\alpha,\alpha}^3(T,T)$.
\ethm
\bpf
Suppose that $\phi_t(x,y,z)=\sum_{i\geq0}\phi_i(x,y,z)t^i$ is the formal $\F[[t]]$-module isomorphism such that $\phi_t\circ\alpha=\alpha\circ\phi_t$ and
$$\sum_{i\geq0}\phi_i\left(\sum_{j\geq0} d_j(x,y,z)t^j\right)t^i =\sum_{i\geq0}d_i'\left(\sum_{k\geq0}\phi_k(x)t^k,\sum_{l\geq0}\phi_l(y)t^l,\sum_{m\geq0}\phi_m(z)t^m\right)t^i.$$
It follows that
$$\sum_{i+j=n}\phi_i(d_j(x,y,z))t^{i+j}=\sum_{i+k+l+m=n}d_i'(\phi_k(x),\phi_l(y),\phi_m(z))t^{i+k+l+m}.$$
In particular,
$$\sum_{i+j=1}\phi_i(d_j(x,y,z))=\sum_{i+k+l+m=1}d_i'(\phi_k(x),\phi_l(y),\phi_m(z)),$$
that is,
$$d_1(x,y,z)+\phi_1([xyz])=[\phi_1(x)yz]+[x\phi_1(y)z]+[xy\phi_1(z)]+d_1'(x,y,z).$$
Then $d_1-d_1'=\delta_{hom}^1\phi_1\in B_{\alpha,\alpha}^3(T,T)$.
\epf
\bdefn
A 3-Hom-cocycle $d_1\in Z_{\alpha,\alpha}^3(T,T)$ is called integrable, if there is a 1-parameter formal deformation $d_t$ of $(T, [\cdot,\cdot,\cdot], \alpha)$ such that $d_t=d_0+d_1t+d_2t^2+\cdots.$
\edefn
\bthm
If $(T, [\cdot,\cdot,\cdot], \alpha)$ is a multiplicative Hom-LTS with $H_{\alpha,\alpha}^5(T,T)=0$, then every 3-Hom-cocycle $d_1\in Z_{\alpha,\alpha}^3(T,T)$ is integrable.
\ethm
\bpf
We proceed by induction on the number of items of $d_t$. First let $d_0=[\cdot,\cdot,\cdot]$. Now suppose that we have already had $d_1,\cdots,d_n\in C_{\alpha,\alpha}^3(T,T)$ satisfying
$$-\delta_{hom}^3d_m=d_1\circ_{\alpha}d_{m-1}+d_2\circ_{\alpha}d_{m-2}+\cdots+d_{m-1}\circ_{\alpha}d_1, \quad\forall m=1,\cdots,n.$$
Set $\tilde{d}=d_1\circ_{\alpha}d_{n}+d_2\circ_{\alpha}d_{n-1}+\cdots+d_{n}\circ_{\alpha}d_1.$ We claim that $\tilde{d}\in Z_{\alpha,\alpha}^5(T,T)$, i.e., $\delta_{hom}^5 \tilde{d}=0$. In fact, for $f,g\in C_{\alpha,\alpha}^3(T,T)$ and $h\in C_{\alpha,\alpha}^5(T,T)$, let
\begin{align*}
&h\bullet_{\alpha}g(x_1,\cdots, x_7)\\
=&\sum_{k=1}^3\sum_{j=2k+1}^7(-1)^{k+1}
h(\alpha(x_1),\cdots,\widehat{x_{2k-1}},\widehat{x_{2k}},\cdots,g(x_{2k-1}x_{2k}x_j),\cdots,\alpha(x_7))
\end{align*}
and
\begin{multline*}
f\bullet_{\alpha}h(x_1,\cdots, x_7)\\\!=\!\!\sum_{k=1}^2\!\sum_{l=k\!+\!1}^3\!\sum_{j=2l\!+\!1}^7\!\!(-\!1)^{k\!+\!l}
f(\alpha^2\!(x_1)\!,\!\cdots\!,\!\widehat{x_{2k\!-\!1}},\!\widehat{x_{2k}},\!\cdots\!,\!\widehat{x_{2l\!-\!1}},\!\widehat{x_{2l}}, \!\cdots\!,\!h(x_{2k\!-\!1}x_{2k}x_{2l\!-\!1}x_{2l}x_j)\!,\!\cdots\!,\!\alpha^2\!(x_7)).
\end{multline*}
Then
\begin{align*}
&\delta_{hom}^5(f\circ_{\alpha}g)(x_1,\cdots, x_7)\\
=&\delta_{hom}^3 f\bullet_{\alpha}g(x_1,\cdots, x_7)+f\bullet_{\alpha}\delta_{hom}^3 g(x_1,\cdots, x_7)\\
 &-\alpha[f(x_1,x_2,x_3)\alpha(x_4) g(x_5,x_6,x_7)]+\alpha[g(x_1,x_2,x_3)\alpha(x_4)f(x_5,x_6,x_7)]\\
 &-\alpha[\alpha(x_3)f(x_1,x_2,x_4) g(x_5,x_6,x_7)]+\alpha[\alpha(x_3)g(x_1,x_2,x_4)f(x_5,x_6,x_7)]\\
 &+\alpha[f(x_1,x_2,x_5)\alpha(x_6) g(x_3,x_4,x_7)]-\alpha[g(x_1,x_2,x_5)\alpha(x_6)f(x_3,x_4,x_7)]\\
 &+\alpha[f(x_1,x_2,x_5) g(x_3,x_4,x_6)\alpha(x_7)]-\alpha[g(x_1,x_2,x_5)f(x_3,x_4,x_6)\alpha(x_7)]\\
 &-\alpha[f(x_1,x_2,x_6)\alpha(x_5) g(x_3,x_4,x_7)]+\alpha[g(x_1,x_2,x_6)\alpha(x_5)f(x_3,x_4,x_7)]\\
 &-\alpha[f(x_1,x_2,x_6)g(x_3,x_4,x_5)\alpha(x_7)]+\alpha[g(x_1,x_2,x_6)f(x_3,x_4,x_5)\alpha(x_7)]\\
 &-\alpha[f(x_3,x_4,x_5)\alpha(x_6) g(x_1,x_2,x_7)]+\alpha[g(x_3,x_4,x_5)\alpha(x_6)f(x_1,x_2,x_7)]\\
 &-\alpha[\alpha(x_5)f(x_3,x_4,x_6) g(x_1,x_2,x_7)]+\alpha[\alpha(x_5)g(x_3,x_4,x_6)f(x_1,x_2,x_7)],
\end{align*}
so
\begin{align*}
\delta_{hom}^5\tilde{d}&=\sum_{i+j=n+1}\delta_{hom}^5(d_i\circ_{\alpha}d_j) =\sum_{i+j=n+1}(\delta_{hom}^3 d_i\bullet_{\alpha}d_j+d_i\bullet_{\alpha}\delta_{hom}^3 d_j)\\
               &=\sum_{i+j+k=n+1}((d_i\circ_{\alpha}d_j)\bullet_{\alpha}d_k+d_i\bullet_{\alpha}(d_j\circ_{\alpha}d_k)).
\end{align*}
Note that
\begin{align*}
&((d_i\circ_{\alpha}d_j)\bullet_{\alpha}d_k+d_i\bullet_{\alpha}(d_j\circ_{\alpha}d_k))(x_1,\cdots, x_7)\\
=&\alpha d_i(d_j(x_1,x_2,x_3),\alpha(x_4), d_k(x_5,x_6,x_7))-\alpha d_i(d_k(x_1,x_2,x_3),\alpha(x_4),d_j(x_5,x_6,x_7))\\
 &+\alpha d_i(\alpha(x_3),d_j(x_1,x_2,x_4), d_k(x_5,x_6,x_7))-\alpha d_i(\alpha(x_3),d_k(x_1,x_2,x_4),d_j(x_5,x_6,x_7))\\
 &-\alpha d_i(d_j(x_1,x_2,x_5),\alpha(x_6), d_k(x_3,x_4,x_7))+\alpha d_i(d_k(x_1,x_2,x_5),\alpha(x_6),d_j(x_3,x_4,x_7))\\
 &-\alpha d_i(d_j(x_1,x_2,x_5), d_k(x_3,x_4,x_6),\alpha(x_7))+\alpha d_i(d_k(x_1,x_2,x_5),d_j(x_3,x_4,x_6),\alpha(x_7))\\
 &+\alpha d_i(d_j(x_1,x_2,x_6),\alpha(x_5), d_k(x_3,x_4,x_7))-\alpha d_i(d_k(x_1,x_2,x_6),\alpha(x_5),d_j(x_3,x_4,x_7))\\
 &+\alpha d_i(d_j(x_1,x_2,x_6),d_k(x_3,x_4,x_5),\alpha(x_7))-\alpha d_i(d_k(x_1,x_2,x_6),d_j(x_3,x_4,x_5),\alpha(x_7))\\
 &+\alpha d_i(d_j(x_3,x_4,x_5),\alpha(x_6), d_k(x_1,x_2,x_7))-\alpha d_i(d_k(x_3,x_4,x_5),\alpha(x_6),d_j(x_1,x_2,x_7))\\
 &+\alpha d_i(\alpha(x_5),d_j(x_3,x_4,x_6), d_k(x_1,x_2,x_7))-\alpha d_i(\alpha(x_5),d_k(x_3,x_4,x_6),d_j(x_1,x_2,x_7)).
\end{align*}
Hence $\delta_{hom}^5 \tilde{d}=0$. Since $H_{\alpha,\alpha}^5(T,T)=0$, it follows that $\tilde{d}\in Z_{\alpha,\alpha}^5(T,T)=B_{\alpha,\alpha}^5(T,T)$. Then there exists $d_{n+1}\in C_{\alpha,\alpha}^3(T,T)$ such that $\tilde{d}=-\delta_{hom}^3 d_{n+1}$. Then we have $d_0,d_1,\cdots,d_n\in C_{\alpha,\alpha}^3(T,T)$ satisfying
$$-\delta_{hom}^3d_m=d_1\circ_{\alpha}d_{m-1}+d_2\circ_{\alpha}d_{m-2}+\cdots+d_{m-1}\circ_{\alpha}d_1, \quad\forall m=1,\cdots,n+1.$$
By induction, we could construct a 1-parameter formal deformation $d_t=d_0+d_1t+d_2t^2+\cdots.$ Therefore, $d_1$ is integrable.
\epf


\begin{thebibliography}{99}
\bibitem{Ammar&Ejbehi&Makhlouf} F. Ammar, Z. Ejbehi and A. Makhlouf, Representations and cohomology of $n$-ary multiplicative Hom-Nambu-Lie algebras. J. Geom. Phys. 61 (2011), no. 10, 1898--1913.
\bibitem{Ammar&Mabrouk&Makhlouf} F. Ammar, S. Mabrouk and A. Makhlouf, Cohomology and deformations of Hom-algebras. J. Lie Theory 21 (2011), no. 4, 813--836.
\bibitem{Benayadi&Makhlouf} S. Benayadi and A. Makhlouf, Hom-Lie algebras with invariant nondegenerate bilinear forms, arXiv: 1009.4226v1 (2010).
\bibitem{Elhamdadi&Makhlouf} M. Elhamdadi and A. Makhlouf, Deformations of Hom-alternative and Hom-Malcev algebras. Algebras Groups Geom. 28 (2011), no. 2, 117--145.
\bibitem{Gerstenhaber1} M. Gerstenhaber, On the deformation of rings and algebras. Ann. of Math. (2) 79 (1964), 59--103.
\bibitem{Gerstenhaber2} M. Gerstenhaber, On the deformation of rings and algebras.II. Ann. of Math. 84 (1966), 1--19.
\bibitem{Gerstenhaber3} M. Gerstenhaber, On the deformation of rings and algebras.III. Ann. of Math. (2) 88 (1968), 1--34.
\bibitem{Gerstenhaber4} M. Gerstenhaber, On the deformation of rings and algebras.IV. Ann. of Math. (2) 99 (1974), 257--276.
\bibitem{Hartwig&Larsson&Silvestrov} J. Hartwig, D. Larsson and S. Silvestrov, Deformations of Lie algebras using $\sigma$-derivations. J. Algebra 295 (2006), no. 2, 314--361.
\bibitem{Hochschild} G. Hochschild, On the cohomology groups of an associative algebra. Ann. of Math. 46 (1945), no. 1, 58--67.
\bibitem{Hu} N. Hu, $Q$-Witt algebras, $q$-Lie algebras, $q$-holomorph structure and representations. Algebra Colloq. 6 (1999), no. 1, 51--70.
\bibitem{Kubo&Taniguchi} F. Kubo and Y.  Taniguchi, A controlling cohomology of the deformation theory of Lie triple systems. J. Algebra 278 (2004), no. 1, 242--250.
\bibitem{Larsson&Silvestrov} D. Larsson and S. Silvestrov, Quasi-hom-Lie algebras, central extensions and 2-cocycle-like identities. J. Algebra 288 (2005), no. 2, 321--344.
\bibitem{Lister} W. Lister, A structure theory of Lie triple systems. Trans. Amer. Math. Soc. 72, (1952), 217--242.
\bibitem{Liu&Chen&Ma} Y. Liu, L. Chen and Y. Ma, Hom-Nijienhuis operators and T*-extensions of hom-Lie superalgebras. Linear Algebra Appl. 439 (2013), no. 7, 2131--2144.
\bibitem{Makhlouf&Silvestrov} A. Makhlouf and S. Silvestrov, Notes on 1-parameter formal deformations of Hom-associative and Hom-Lie algebras. Forum Math. 22 (2010), no. 4, 715--739.
\bibitem{Sheng} Y. Sheng, Representations of hom-Lie algebras.  Algebr. Represent. Theory 15 (2012), no. 6, 1081--1098.
\bibitem{Sheng&Chen} Y. Sheng and D. Chen, Hom-Lie 2-algebras. J. Algebra 376 (2013), 174--195.
\bibitem{Yamaguti} K. Yamaguti, On the cohomology space of Lie triple system. Kumamoto J. Sci. Ser. A 5 (1960),  44--52.
\bibitem{Yau1} D. Yau, Hom-algebras and homology. J. Lie Theory 19 (2009), no. 2, 409--421.
\bibitem{Yau2} D. Yau, On $n$-ary Hom-Nambu and Hom-Nambu-Lie algebras. J. Geom. Phys. 62 (2012), no. 2, 506--522.
\end{thebibliography}
\end{document}